\documentclass[fleqn]{mat01}
\usepackage{times,mathtimy,amssymb,latexsym}
\begin{document}

\setcounter{page}{167}
\firstpage{167}

\newcommand{\sign}{\mbox{sign}}
\newcommand{\Discr}{\mbox{Discr}}
\newcommand{\transrow}[1]{\stackrel{#1}{\longrightarrow}}

\newtheorem{theor}{\bf Theorem}
\newtheorem{lem}{\it Lemma}
\newtheorem{application}{\it Application}
\def\definit{\trivlist \item[\hskip \labelsep{\rm DEFINITION}]}

\newtheorem{theore}{Theorem}
\renewcommand\thetheore{1({\it \alph{theore}})}
\newtheorem{pot}[theore]{\it Proof of Theorem}

\newtheorem{theoree}{Theorem}
\renewcommand\thetheoree{2({\it \roman{theoree}})}
\newtheorem{pott}[theoree]{\it Proof of Theorem}

\title{Some BMO estimates for vector-valued multilinear
singular integral operators}

\markboth{Liu Lanzhe}{Some BMO estimates}

\author{LIU LANZHE}

\address{College of Mathematics and Computer, Changsha University of
Science and~Technology, Changsha 410077, People's Republic of China\\
\noindent E-mail: lanzheliu@263.net}

\volume{115}

\mon{May}

\parts{2}

\pubyear{2005}

\Date{MS received 4 April 2004; revised 19 December 2004}

\begin{abstract}
In this paper, we prove some BMO end-point estimates for some
vector-valued multilinear operators related to certain singular integral
operators.
\end{abstract}


\maketitle

\section{Introduction and notations}

Let $b\in \hbox{BMO}(R^n)$ and $T$ be the Calder\'on--Zygmund singular
integral operator. The commutator $[b, T]$ generated by $b$ and $T$ is
defined as $[b, T](f)(x)=b(x)T(f)(x)-T(bf)(x)$. By using a classical
result of Coifman {\it et~al} \cite{8}, we know that the commutator
$[b,T]$ is bounded on $L^p(R^n)$ for $1<p<\infty$. Chanillo \cite{1}
proves a similar result when $T$ is replaced by the fractional integral
operator. However, it was observed that the commutator is not bounded,
in general, from $H^p(R^n)$ to $L^p(R^n)$ for $0<p\leq 1$
\cite{13,14,15}. In \cite{11}, the boundedness properties of the
commutator for the extreme values of $p$ are obtained. Also, in
\cite{2}, Chanillo studies some commutators generated by a very general
class of pseudo-differential operators and proves the boundedness on
$L^p(R^n)\ (1<p<\infty)$ for the commutators, and note that the conditions
on the kernel of the singular integral operator arise from a
pseudo-differential operator. As the development of singular integral
operators and their commutators, multilinear singular integral operators
have been well-studied. It is known that multilinear operator, as a
non-trivial extension of the commutator, is of great interest in
harmonic analysis and has been widely studied by many authors 
\cite{3,4,5,6,7}. In \cite{9}, the weighted $L^p(p>1)$-boundedness of
the multilinear operator related to some singular integral operators is
obtained and in \cite{3}, the weak ($H^1$, $L^1$)-boundedness of the
multilinear operator related to some singular integral operators is
obtained. The main purpose of this paper is to establish the BMO
end-point estimates for some vector-valued multilinear operators related
to certain singular integral operators. 

First, let us introduce some notations \cite{10,16}. Throughout this
paper, $Q=Q(x, r)$ will denote a cube of $R^n$ with sides parallel to
the axes and centered at $x$ and having side length. For a locally
integrable function $f$ and non-negative weight function $w$, let 
\hbox{$w(Q)=\int_Qw(x){\rm d}x$}, $f_{w, Q}=w(Q)^{-1}\int_Qf(x)w(x){\rm d}x$ and
$f^{\#}(x) = \sup_{x\in Q}w(Q)^{-1}\int_Q$ $|f(y)-f_{w, Q}|w(x){\rm d}y$.
$f$ is said to belong to $\hbox{BMO}(w)$ if $f^{\#}\in L^\infty(w)$ and
define $\|f\|_{{\rm BMO}(w)}=\|f^{\#}\|_{L^\infty(w)}$. We denote
$\hbox{BMO}(w)=\hbox{BMO}(R^n)$ and $\|f\|_{\rm
BMO}=\|f^{\#}\|_{L^\infty}$ if $w=1$. It is well-known that \cite{12}
\begin{equation*}
\|f\|_{{\rm BMO}(w)} \approx \sup\limits_Q \inf\limits_{c\in C}
w(Q)^{-1} \int_Q|f(x)-c|w(x){\rm d}x.
\end{equation*}
We also define the weighted central BMO space by $\hbox{CMO}(w)$, which
is the space of those functions $f\in L_{\rm loc}(R^n)$ such that
\begin{equation*}
\|f\|_{{\rm CMO}(w)} = \sup\limits_{d>1} w(Q(0,d))^{-1} \int_Q
|f(y)-f_{w,Q}|w(y){\rm d}y < \infty.
\end{equation*}
We denote the Muckenhoupt weights by $A_p$ for $1\leq p<\infty$
\cite{10}, that is,
\begin{align*}
A_p &= \left\{0<w\in L^1_{\rm loc}(R^n)\hbox{:}\
\sup_Q\left(\frac{1}{|Q|}\int_Q w(x){\rm d}x\right) \right.\\[.3pc]
&\quad\, \left. \times \left(\frac{1}{|Q|}\int_Q w(x)^{-1/(p-1)}{\rm
d}x\right)^{p-1}<\infty\right\},\quad 1<p<\infty,\\[.3pc]
A_1 &= \left\{0<w\in L^1_{\rm loc}(R^n)\hbox{:}\ \sup_{x\in
Q}\frac{w(Q)}{|Q|}\leq Cw(x), \hbox{a.e.} \right\}
\end{align*}
and
\begin{equation*}
A_{\infty} = \bigcup_{1\leq p<\infty}A_p.
\end{equation*}

\vspace{.1pc}
\begin{definit}$\left.\right.$
{\rm 
\begin{enumerate}
\renewcommand\labelenumi{(\arabic{enumi})}
\item Let $0<\delta<n$ and $1<p<n/\delta$. We shall call
$B_p^{\delta}(R^n)$ the space of those functions $f$ on $R^n$ such that
\begin{equation*}
\hskip -1.25pc \|f\|_{B_p^{\delta}}=\sup\limits_{r>1}r^{-n(1/p-\delta/n)}\|f\chi_{Q(0,
r)}\|_{L^p}<\infty.
\end{equation*}

\item Let $1<p<\infty$ and $w$ be a non-negative weight function on
$R^n$. We shall call $B_p(w)$ the space of that function $f$ on $R^n$
such that
\begin{equation*}
\hskip -1.25pc \|f\|_{B_p(w)}=\sup\limits_{r>1}[w(Q(0, r))]^{-1/p}\|f\chi_{Q(0,
r)}\|_{L^p(w)}<\infty.
\end{equation*}
\end{enumerate}}\vspace{-1pc}
\end{definit}

\section{Theorems}

In this paper, we will study a class of vector-valued multilinear
operators related to some singular integral operators, whose definitions
are the following.

Fix $\varepsilon>0$ and $\delta\geq 0$. Let $T\hbox{:}\ S\to S'$ be a
linear operator and there exists a locally integrable function $K(x,y)$
on $R^n \times R^n\setminus \{(x,y)\in R^n\times R^n\hbox{:}\ x=y\}$
such that
\begin{equation*}
T_\delta(g)(x) = \int_{R^n}K(x, y)g(y){\rm d}y
\end{equation*}
for every bounded and compactly supported function $g$, where $K$
satisfies:
\begin{equation*}
|K(x,y)|\leq C|x-y|^{-n+\delta}
\end{equation*}
and
\begin{equation*}
|K(y,x)-K(z,x)|+|K(x,y)-K(x,z)|\leq
C|y-z|^\varepsilon|x-z|^{-n-\varepsilon+\delta}
\end{equation*}
if $2|y-z|\leq |x-z|$. Let $m_j$ be the positive integers
($j=1,\ldots,l$), $m_1+\cdots+m_l=m$ and $A_j$ be the functions on
$R^n\,(j=1,\ldots,l)$. For $1<r<\infty$, the vector-valued multilinear
operator associated with $T$ is defined as
\begin{equation*}
|T_\delta^A(f)(x)|_r=\left(\sum_{i=1}^\infty
|T_\delta^A(f_i)(x)|^r\right)^{1/r},
\end{equation*}
where
\begin{equation*}
T_\delta^A(f_i)(x)=\int_{R^n}\frac{\prod_{j=1}^lR_{m_j+1}(A_j; x,
y)}{|x-y|^m}K(x,y)f_i(y){\rm d}y
\end{equation*}
and
\begin{equation*}
R_{m_j+1}(A_j; x,y)=A_j(x)-\sum_{|\alpha|\leq
m_j}\frac{1}{\alpha!}D^\alpha A_j(y)(x-y)^\alpha.
\end{equation*}
Set
\begin{align*}
|T_\delta(f)(x)|_r=\left(\sum_{i=1}^\infty |T(f_i)(x)|^r\right)^{1/r}\quad
\hbox{and}\quad |f|_r=\left(\sum_{i=1}^\infty |f_i(x)|^r\right)^{1/r}.
\end{align*}
We write $T_\delta=T$, $|T_\delta|_r=|T|_r$ and $|T_\delta^A|_r=|T^A|_r$
if $\delta=0$.

Note that when $m=0$, $T_\delta^A$ is just the multilinear commutators
of $T_\delta$ and $A$ \cite{13,14,15}. In this paper, we will prove the
BMO estimates for the vector-valued multilinear operators
$|T_\delta^A|_r$ and $|T^A|_r$.

Now we state our results as follows.

\begin{theor}[\!] Let $1<r<\infty, 0<\delta<n, 1<p<n/\delta$ and
$D^\alpha A_j\in {\rm BMO}(R^n)$ for all $\alpha$ with $|\alpha|=m_j$
and $j=1,\ldots,l$. Suppose that $|T_\delta|_r$ maps $L^s(R^n)$
continuously into $L^t(R^n)$ for any $s,t\in (1,+\infty]$ with
$1<s<n/\delta$ and $1/t=1/s-\delta/n$. Then

\begin{enumerate}
\renewcommand\labelenumi{\rm (\alph{enumi})}
\item $|T_\delta^A|_r$ maps $L^{n/\delta}(R^n)$ continuously into ${\rm
BMO}(R^n)${\rm ,} that is
\begin{equation*}
\hskip -1.25pc \||T_\delta^A(f)|_r\|_{\rm BMO}\leq C\||f|_r\|_{L^{n/\delta}}.
\end{equation*}

\item $|T_\delta^A|_r$ maps $B_p^{\delta}(R^n)$ continuously into
${\rm CMO}(R^n)${\rm ,} that is
\begin{equation*}
\hskip -1.25pc \||T_\delta^A(f)|_r\|_{\rm CMO}\leq C\||f|_r\|_{B_p^{\delta}}.
\end{equation*}
\end{enumerate}
\end{theor}

\begin{theor}[\!] Let $1<r<\infty, 1<p<\infty$ and $D^\alpha A_j\in {\rm
BMO}(R^n)$ for all $\alpha$ with $|\alpha|=m_j$ and $j=1,\ldots,l$.

\begin{enumerate}
\renewcommand\labelenumi{\rm (\roman{enumi})}
\item If $w \in A_\infty$ and that $|T|_r$ is bounded on $L^s(w)$ for
any $1<s\leq\infty$ and $w\in A_\infty${\rm ,} then $|T^A|_r$ maps
$L^\infty(w)$ continuously into ${\rm BMO}(w)${\rm ,} that is{\rm ,}
\begin{equation*}
\hskip -1.25pc \||T^A(f)|_r\|_{{\rm BMO}(w)}\leq C\||f|_r\|_{L^\infty(w)};
\end{equation*}

\item If $w \in A_1$ and that $|T|_r$ is bounded on $L^s(w)$ for any
$1<s\leq\infty$ and $w\in A_1${\rm ,} then $|T^A|_r$ maps $B_p(w)$
continuously into ${\rm CMO}(w)${\rm ,} that is{\rm ,}
\begin{equation*}
\hskip -1.25pc \||T^A(f)|_r\|_{{\rm CMO}(w)}\leq C\||f|_r\|_{B_p(w)}.
\end{equation*}
\end{enumerate}
\end{theor}

\section{Proofs of theorems}

To prove the theorems, we need the following lemmas.

\begin{lem}\hskip -.3pc {\rm \cite{6}.}\ \ Let $A$ be a function on $R^n$ and
$D^\alpha A\in L^q(R^n)$ for all $\alpha$ with $|\alpha|=m$ and some
$q>n$. Then
\begin{equation*}
|R_m(A;x,y)|\leq C|x-y|^m\sum_{|\alpha|=m}\left(\frac{1}{|\tilde
Q(x,y)|}\int_{\tilde Q(x,y)}|D^\alpha A(z)|^q {\rm d}z\right)^{1/q},
\end{equation*}

$\left.\right.$\vspace{-1.5pc}

\noindent where $\tilde Q(x,y)$ is the cube centered at $x$ and having side length
$5\sqrt{n}|x-y|$.
\end{lem}

\begin{lem}
Let $w \in A_\infty${\rm ,} then ${\rm BMO}(w)={\rm BMO}(R^n)$.
\end{lem}

The proof of the lemma follows from \cite{12} and the John--Nirenberg
Lemma for BMO \cite{10}.

\begin{pot}{\rm
It is only to prove that there exists a constant $C_Q$ such that
\begin{equation*}
\frac{1}{|Q|}\int_Q||T_\delta^A(f)(x)|_r-C_Q|{\rm d}x \leq
C\||f|_r\|_{L^{n/\delta}}
\end{equation*}
holds for any cube $Q$. Without loss of generality, we may assume $l=2$.
Fix a cube $Q=Q(x_0, d)$. Let $\tilde Q=5\sqrt{n}Q$ and $\tilde
A_j(x)=A_j(x)-\sum_{|\alpha|=m_j}\frac{1}{\alpha!}(D^\alpha
A_j)_{\tilde Q}x^\alpha$, then $R_{m_j}(A_j;x,y) =R_{m_j}(\tilde
A_j;x,y)$ and $D^\alpha\tilde A_j=D^\alpha A_j-(D^\alpha A_j)_{\tilde
Q}$ for $|\alpha|=m_j$. We split $f=g+h=\{g_i\}+\{h_i\}$ for
$g_i=f_i\chi_{\tilde Q}$ and $h_i=f_i\chi_{R^n\setminus\tilde Q}$. Write
\begin{align*}
\hskip -4pc T_\delta^A(f_i)(x)
&= \int_{R^n}\frac{\prod_{j=1}^2R_{m_j+1}(\tilde A_j; x, y)}{|x-y|^m}K(x,y)f_i(y){\rm d}y\\[.5pc]
\hskip -4pc &= \int_{R^n} \frac{\prod_{j=1}^2R_{m_j+1}(\tilde A_j; x, y)}{|x-y|^m}K(x, y)h_i(y){\rm d}y  \\[.5pc]
\hskip -4pc &\quad\, +\int_{R^n}\frac{\prod_{j=1}^2R_{m_j}(\tilde A_j; x, y)}{|x-y|^m}K(x, y)g_i(y){\rm d}y\\[.5pc]
\hskip -4pc &\quad\, -\sum_{|\alpha_1|=m_1}\frac{1}{\alpha_1!}\int_{R^n}\frac{R_{m_2}(\tilde A_2; x, y)(x-y)^{\alpha_1}}{|x-y|^m}D^{\alpha_1}\tilde A_1(y)K(x,y)g_i(y){\rm d}y
\end{align*}
\begin{align*}
\hskip -4pc \phantom{T_\delta^A(f_i)(x)} &\quad\, -\sum_{|\alpha_2|=m_2}\frac{1}{\alpha_2!}\int_{R^n}\frac{R_{m_1}(\tilde A_1; x, y)(x-y)^{\alpha_2}}{|x-y|^m}D^{\alpha_2}\tilde A_2(y)K(x,y)g_i(y){\rm d}y\\[.5pc]
\hskip -4pc \phantom{T_\delta^A(f_i)(x)} &\quad\, +\sum_{|\alpha_1|=m_1, |\alpha_2|=m_2}\frac{1}{\alpha_1!\alpha_2!}\\[.5pc]
\hskip -4pc \phantom{T_\delta^A(f_i)(x)} &\quad\, \times \int_{R^n}\frac{(x-y)^{\alpha_1+\alpha_2}
D^{\alpha_1}\tilde A_1(y)D^{\alpha_2}\tilde A_2(y)}{|x-y|^m}K(x,y)g_i(y){\rm d}y,
\end{align*}
then, by the Minkowski's inequality,
\begin{align*}
&\frac{1}{|Q|}\int_Q ||T_\delta^A(f)(x)|_r-|T_\delta^{\tilde A}(h)(x_0)|_r|{\rm d}x\\[.5pc]
&\quad\,\leq\frac{1}{|Q|}\int_Q\left(\sum_{i=1}^\infty |T_A(f_i)(x)-T_{\tilde A}(h_i)(x_0)|^r\right)^{1/r}{\rm d}x \\[.5pc]
&\quad\,\leq \frac{1}{|Q|}\int_Q\left(\sum_{i=1}^\infty\left|\int_{R^n}\frac{\prod_{j=1}^2R_{m_j}(\tilde A_j; x, y)}{|x-y|^m}K(x, y)g_i(y){\rm d}y\right|^r\right)^{1/r}{\rm d}x \\[.5pc]
&\qquad\,+\frac{1}{|Q|}\int_Q\left(\sum_{i=1}^\infty\left|\sum_{|\alpha_1|=m_1}\frac{1}{\alpha_1!} \right.\right.\\[.5pc]
&\qquad\, \left.\left. \times \int_{R^n}\frac{R_{m_2}(\tilde A_2; x, y)(x-y)^{\alpha_1}}{|x-y|^m}D^{\alpha_1}\tilde A_1(y)K(x,y)g_i(y){\rm d}y\right|^r\right)^{1/r}{\rm d}x \\[.5pc]
&\qquad\,+\frac{1}{|Q|}\int_Q\left(\sum_{i=1}^\infty\left|\sum_{|\alpha_2|=m_2}\frac{1}{\alpha_2!}\right.\right.\\[.5pc]
&\qquad\, \left.\left. \times \int_{R^n}\frac{R_{m_1}(\tilde A_1; x, y)(x-y)^{\alpha_2}}{|x-y|^m}D^{\alpha_2}\tilde A_2(y)K(x,y)g_i(y){\rm d}y\right|^r\right)^{1/r}{\rm d}x \\[.5pc]
&\qquad\,+\frac{1}{|Q|}\int_Q\left(\sum_{i=1}^\infty\left|\sum_{|\alpha_1|=m_1,\ |\alpha_2|=m_2}\frac{1}{\alpha_1!\alpha_2!} \right.\right.\\[.5pc]
&\qquad\, \left.\left. \times \int_{R^n}\frac{(x-y)^{\alpha_1+\alpha_2}D^{\alpha_1}\tilde A_1(y)D^{\alpha_2}\tilde A_2(y)}{|x-y|^m}K(x,y)g_i(y){\rm d}y\right|^r\right)^{1/r}{\rm d}x \\[.5pc]
&\qquad\,+\frac{1}{|Q|}\int_Q\left(\sum_{i=1}^\infty \left|T_\delta^{\tilde A}(h_i)(x)-T_\delta^{\tilde A}(h_i)(x_0)\right|^r\right)^{1/r}{\rm d}x \\[.5pc]
&\quad\,:= {\rm I}_1 + {\rm I}_2 + {\rm I}_3 + {\rm I}_4 + {\rm I}_5.
\end{align*}
Now, let us estimate $\hbox{I}_1, \hbox{I}_2, \hbox{I}_3, \hbox{I}_4$
and $\hbox{I}_5$, respectively. First, for $x\in Q$ and $y\in \tilde Q$,
by Lemma~1, we get
\begin{equation*}
R_{m_j}(\tilde A_j; x, y)\leq C|x-y|^{m_j}
\sum_{|\alpha_j|=m_j}\|D^{\alpha_j}A_j\|_{\rm BMO}.
\end{equation*}
Thus, by the $(L^{n/\delta}, L^\infty)$-boundedness of $|T_\delta|_r$,
we get
\begin{align*}
{\rm I}_1 &\leq C\prod_{j=1}^2\left(\sum_{|\alpha_j|=m_j}\|D^{\alpha_j}A_j\|_{\rm BMO}\right)\frac{1}{|Q|}\int_Q |T_\delta(g)(x)|_r{\rm d}x \\[.4pc]
&\leq C\prod_{j=1}^2\left(\sum_{|\alpha_j|=m_j}\|D^{\alpha_j}A_j\|_{\rm BMO}\right)\||T_\delta(g)|_r\|_{L^\infty} \\[.4pc]
&\leq C\prod_{j=1}^2\left(\sum_{|\alpha_j|=m_j}\|D^{\alpha_j}A_j\|_{\rm BMO}\right)\||f|_r\|_{L^{n/\delta}}.
\end{align*}
For ${\rm I}_2$, by the $(L^p,L^q)$-boundedness of $T_\delta$ for
$1/q=1/p-\delta/n$, $n/\delta>p>1$ and the H\"older's inequality, we get
\begin{align*}
{\rm I}_2 &\leq C\sum_{|\alpha_2|=m_2}\|D^{\alpha_2}A_2\|_{\rm BMO}\sum_{|\alpha_1|=m_1}\frac{1}{|Q|}\int_Q |T_\delta(D^{\beta_1}\tilde A_1g)(x)|_r{\rm d}x \\[.4pc]
&\leq C\sum_{|\alpha_2|=m_2}\|D^{\alpha_2}A_2\|_{\rm BMO}\sum_{|\alpha_1|=m_1}\left(\frac{1}{|Q|}\int_{R^n}|T_\delta(D^{\alpha_1}\tilde A_1g)(x)|_r^q{\rm d}x\right)^{1/q}\\[.4pc]
&\leq C\sum_{|\alpha_2|=m_2}\|D^{\alpha_2}A_2\|_{\rm BMO}\sum_{|\alpha_1|=m_1}|Q|^{-1/q}\left(\int_{R^n}|D^{\alpha_1}\tilde A_1(x)g(x)|_r^p{\rm d}x\right)^{1/p} \\[.4pc]
&\leq C\sum_{|\alpha_2|=m_2}\|D^{\alpha_2}A_2\|_{\rm BMO}\\[.4pc]
&\quad\, \times \sum_{|\alpha_1|=m_1}\left(\frac{1}{|Q|}\int_{\tilde Q}|D^{\alpha_1}A_1(x)-(D^{\alpha_1}A_1)_{\tilde Q}|^q{\rm d}x\right)^{1/q}\||f|_r\|_{L^{n/\delta}}  \\[.4pc]
&\leq C\prod_{j=1}^2\left(\sum_{|\alpha_j|=m_j}\|D^{\alpha_j}A_j\|_{\rm BMO}\right)\||f|_r\|_{L^{n/\delta}}.
\end{align*}
For ${\rm I}_3$, similar to the proof of ${\rm I}_2$, we get
\begin{equation*}
{\rm I}_3 \leq C\prod_{j=1}^2\left(\sum_{|\alpha_j|=m_j}\|D^{\alpha_j}A_j\|_{\rm
BMO}\right)\||f|_r\|_{L^{n/\delta}}.
\end{equation*}
Similarly, for ${\rm I}_4$, choose $1<p<n/\delta$ and $q, t_1,t_2>1$
such that $1/q=1/p-\delta/n$ and $1/t_1+1/t_2+p\delta/n=1$. We obtain,
by the H\"older's inequality,
\begin{align*}
{\rm I}_4 &\leq C\sum_{|\alpha_1|=m_1, |\alpha_2|=m_2}\frac{1}{|Q|}\int_Q |T_\delta(D^{\alpha_1}\tilde A_1D^{\alpha_2}\tilde A_2g)(x)|_r{\rm d}x \\[.4pc]
&\leq C\sum_{|\alpha_1|=m_1, |\alpha_2|=m_2}\left(\frac{1}{|Q|}\int_{R^n}|T_\delta(D^{\alpha_1}\tilde A_1D^{\alpha_2}\tilde A_2g)(x)|_r^q{\rm d}x\right)^{1/q} \\[.4pc]
&\leq C\sum_{|\alpha_1|=m_1, |\alpha_2|=m_2}|Q|^{-1/q}\left(\int_{R^n}|D^{\alpha_1}\tilde A_1(x)D^{\alpha_2}\tilde A_2(x)g(x)|_r^p{\rm d}x\right)^{1/p} \\[.4pc]
&\leq C\sum_{|\alpha_1|=m_1, |\alpha_2|=m_2}\left(\frac{1}{|Q|}\int_{\tilde Q}|D^{\alpha_1}\tilde A_1(x)|^{pt_1}{\rm d}x\right)^{1/pt_1}\\[.4pc]
&\quad\, \times \left(\frac{1}{|Q|}\int_{\tilde Q}|D^{\alpha_2}\tilde A_2(x)|^{pt_2}{\rm d}x\right)^{1/pt_2}\||f|_r\|_{L^{n/\delta}} \\[.4pc]
&\leq C\prod_{j=1}^2\left(\sum_{|\alpha_j|=m_j}\|D^{\alpha_j}A_j\|_{\rm BMO}\right)\||f|_r\|_{L^{n/\delta}}.
\end{align*}
For ${\rm I}_5$, we write
\begin{align*}
\hskip -4pc &T_\delta^{\tilde A}(h_i)(x)-T_\delta^{\tilde A}(h_i)(x_0)\\[.3pc]
\hskip -4pc &\quad\, =\int_{R^n}\left(\frac{K(x,y)}{|x-y|^m}-\frac{K(x_0,y)}{|x_0-y|^m}\right)\prod_{j=1}^2R_{m_j}(\tilde A_j; x, y)h_i(y){\rm d}y \\[.4pc]
\hskip -4pc &\qquad\, +\int_{R^n}\left(R_{m_1}(\tilde A_1; x, y)-R_{m_1}(\tilde A_1; x_0, y)\right)\frac{R_{m_2}(\tilde A_2; x, y)}{|x_0-y|^m}K(x_0, y)h_i(y){\rm d}y \\[.4pc]
\hskip -4pc &\qquad\, +\int_{R^n}\left(R_{m_2}(\tilde A_2; x, y)-R_{m_2}(\tilde A_2; x_0, y)\right)\frac{R_{m_1}(\tilde A_1; x_0, y)}{|x_0-y|^m}K(x_0, y)h_i(y){\rm d}y\\[.4pc]
\hskip -4pc &\qquad\, -\sum_{|\alpha_1|=m_1}\frac{1}{\alpha_1!}\int_{R^n}\left[\frac{R_{m_2}(\tilde A_2; x, y)(x-y)^{\alpha_1}}{|x-y|^m}K(x,y) \right.\\[.4pc]
\hskip -4pc &\qquad\ \!\left. -\frac{R_{m_2}(\tilde A_2; x_0, y)(x_0-y)^{\alpha_1}}{|x_0-y|^m}K(x_0,y)\right]D^{\alpha_1}\tilde A_1(y)h_i(y){\rm d}y\\[.4pc]
\hskip -4pc &\qquad\, -\sum_{|\alpha_2|=m_2}\frac{1}{\alpha_2!}\int_{R^n}\left[\frac{R_{m_1}(\tilde A_1; x, y)(x-y)^{\alpha_2}}{|x-y|^m}K(x,y) \right.\\[.4pc]
\hskip -4pc &\qquad\ \!\left. -\frac{R_{m_1}(\tilde A_1; x_0, y)(x_0-y)^{\alpha_2}}{|x_0-y|^m}K(x_0,y)\right]D^{\alpha_2}\tilde A_2(y)h_i(y){\rm d}y
\end{align*}
\begin{align*}
\hskip -4pc &\qquad\, +\sum_{|\alpha_1|=m_1,\ |\alpha_2|=m_2}\frac{1}{\alpha_1!\alpha_2!}\int_{R^n}\left[\frac{(x-y)^{\alpha_1+\alpha_2}}{|x-y|^m}K(x,y) \right.\\[.4pc]
\hskip -4pc &\qquad\, \left. -\frac{(x_0-y)^{\alpha_1+\alpha_2}}{|x_0-y|^m}K(x_0,y)\right]D^{\alpha_1}\tilde A_1(y)D^{\alpha_2}\tilde A_2(y)h_i(y){\rm d}y  \\[.4pc]
\hskip -4pc &\quad\, = {\rm I}_5^{(1)} + {\rm I}_5^{(2)}+{\rm I}_5^{(3)}+{\rm I}_5^{(4)}+{\rm I}_5^{(5)}+{\rm I}_5^{(6)}.
\end{align*}
By Lemma~1 and the following inequality \cite{16}
\begin{equation*}
|b_{Q_1}-b_{Q_2}|\leq C\log(|Q_2|/|Q_1|)\|b\|_{\rm BMO}\quad \hbox{for}\
Q_1 \subset Q_2,
\end{equation*}
we know that, for $x\in Q$ and $y\in 2^{k+1}\tilde Q\setminus 2^k\tilde Q$,
\begin{align*}
\hskip -4pc |R_{m_j}(\tilde A_j;x,y)|&\leq C|x-y|^{m_j}\sum_{|\alpha|=m_j}(\|D^\alpha A_j\|_{\rm BMO}+|(D^\alpha A_j)_{\tilde Q(x,y)}-(D^\alpha A_j)_{\tilde Q}|) \\[.3pc]
&\leq C k|x-y|^{m_j}\sum_{|\alpha|=m_j}\|D^\alpha A_j\|_{\rm BMO}.
\end{align*}
Note that $|x-y|\sim |x_0-y|$ for $x\in Q$ and $y\in R^n\setminus\tilde
Q$, and we obtain, by the conditions on $K$,
\begin{align*}
\hskip -4pc |{\rm I}_5^{(1)}|&\leq C\int_{R^n}\left(\frac{|x-x_0|}{|x_0-y|^{m+n+1-\delta}}+\frac{|x-x_0|^\varepsilon}
{|x_0-y|^{m+n+\varepsilon-\delta}}\right)\prod_{j=1}^2|R_{m_j}(\tilde A_j; x, y)||h_i(y)|{\rm d}y \\[.3pc]
\hskip -4pc &\leq C\prod_{j=1}^2\left(\sum_{|\alpha_j|=m_j}\|D^{\alpha_j}A_j\|_{\rm
BMO}\right)\\[.3pc]
\hskip -4pc &\quad\, \times \sum_{k=0}^\infty\int_{2^{k+1}\tilde Q\setminus2^k\tilde
Q}k^2\left(\frac{|x-x_0|}
{|x_0-y|^{n+1-\delta}}+\frac{|x-x_0|^\varepsilon}{|x_0-y|^{n+\varepsilon-
\delta}}\right)|f_i(y)|{\rm d}y.
\end{align*}
Thus, by the Minkowski's inequality,
\begin{align*}
\hskip -4pc \left(\sum_{i=1}^\infty|{\rm I}_5^{(1)}|^r\right)^{1/r}
&\leq C\prod_{j=1}^2\left(\sum_{|\alpha_j|=m_j}\|D^{\alpha_j}A_j\|_{\rm BMO}\right)\\[.3pc]
\hskip -4pc &\quad\, \times \sum_{k=0}^\infty\int_{2^{k+1}\tilde Q\setminus2^k\tilde Q}k^2\left(\frac{|x-x_0|}{|x_0-y|^{n+1-\delta}}
+\frac{|x-x_0|^\varepsilon}{|x_0-y|^{n+\varepsilon-\delta}}\right)|f(y)|_r{\rm d}y  \\[.3pc]
\hskip -4pc &\leq C\prod_{j=1}^2\left(\sum_{|\alpha_j|=m_j}\|D^{\alpha_j}A_j\|_{\rm BMO}\right)
\sum_{k=1}^\infty k^2(2^{-k}+2^{-\varepsilon k})\||f|_r\|_{L^{n/\delta}} \\[.3pc]
\hskip -4pc &\leq C\prod_{j=1}^2\left(\sum_{|\alpha_j|=m_j}\|D^{\alpha_j}A_j\|_{\rm BMO}\right)\||f|_r\|_{L^{n/\delta}}.
\end{align*}
For ${\rm I}_5^{(2)}$, by the formula \cite{6}:
\begin{equation*}
R_{m_j}(\tilde A_j; x, y)-R_{m_j}(\tilde A_j; x_0,
y)=\sum_{|\beta|<m}\frac{1}{\beta!}R_{m-|\beta|}(D^\beta\tilde A_j; x,
x_0)(x-y)^\beta
\end{equation*}

$\left.\right.$\vspace{-1.7pc}

\noindent and Lemma~1, we have
\begin{align*}
&|R_{m_j}(\tilde A_j; x, y)-R_{m_j}(\tilde A_j; x_0, y)|\\[.3pc]
&\quad\, \leq C\sum_{|\beta|<m_j}\sum_{|\alpha|=m_j}|x-x_0|^{m_j-|\beta|}|x-
y|^{|\beta|}\|D^\alpha A_j\|_{\rm BMO}.
\end{align*}
Thus
\begin{align*}
\left(\sum_{i=1}^\infty|{\rm I}_5^{(2)}|^r\right)^{1/r}&\leq C\prod_{j=1}^2\left(\sum_{|\alpha|=m_j}
\|D^\alpha A_j\|_{\rm BMO}\right)\\[.3pc]
&\quad\, \times \sum_{k=0}^\infty\int_{2^{k+1}\tilde Q\setminus2^k\tilde Q}k\frac{|x-x_0|}{|x_0-y|^{n+1-\delta}}|f(y)|_r{\rm d}y  \\[.3pc]
&\leq C\prod_{j=1}^2\left(\sum_{|\alpha_j|=m_j}\|D^{\alpha_j}A_j\|_{\rm BMO}\right)\||f|_r\|_{L^{n/\delta}}.
\end{align*}
Similarly,
\begin{equation*}
\left(\sum_{i=1}^\infty|{\rm I}_5^{(3)}|^r\right)^{1/r}\leq
C\prod_{j=1}^2\left(\sum_{|\alpha_j|=m_j}\|D^{\alpha_j}A_j\|_{\rm
BMO}\right)\||f|_r\|_{L^{n/\delta}}.
\end{equation*}
For ${\rm I}_5^{(4)}$, taking $t>1$ such that $1/t+\delta/n=1$, then
\begin{align*}
\hskip -4pc \left(\sum_{i=1}^\infty|{\rm I}_5^{(4)}|^r\right)^{1/r}
&\leq C\sum_{|\alpha_1|=m_1}\int_{R^n}\left|\frac{(x-y)^{\alpha_1}K(x,y)}{|x-y|^m}-\frac{(x_0-y)^{\alpha_1}K(x_0,y)}
{|x_0-y|^m}\right|\\[.3pc]
\hskip -4pc &\quad\,\times |R_{m_2}(\tilde A_2; x, y)||D^{\alpha_1}\tilde A_1(y)||h(y)|_r{\rm d}y  \\[.3pc]
\hskip -4pc &\quad\,+C\sum_{|\alpha_1|=m_1}\int_{R^n}|R_{m_2}(\tilde A_2; x, y)-R_{m_2}(\tilde A_2; x_0, y)|\\[.2pc]
\hskip -4pc &\quad\, \times\frac{|(x_0-y)^{\alpha_1}K(x_0,y)|}
{|x_0-y|^m}|D^{\alpha_1}\tilde A_1(y)||h(y)|_r{\rm d}y\\[.3pc]
\hskip -4pc &\leq C\sum_{|\alpha_2|=m_2}\|D^{\alpha_2}A_2\|_{\rm BMO}\sum_{|\alpha_1|=m_1}\sum_{k=1}^\infty k(2^{-k}+2^{-\varepsilon k})\\[.3pc]
\hskip -4pc &\quad\,\times \left(\frac{1}{|2^k\tilde Q|}\int_{2^k\tilde Q}|D^{\alpha_1}\tilde A_1(y)|^t{\rm d}y\right)^{1/t}\||f|_r\|_{L^{n/\delta}} \\[.3pc]
\hskip -4pc &\leq C\prod_{j=1}^2\left(\sum_{|\alpha_j|=m_j}\|D^{\alpha_j}A_j\|_{\rm BMO}\right)\||f|_r\|_{L^{n/\delta}}.
\end{align*}
Similarly,
\begin{equation*}
\left(\sum_{i=1}^\infty|{\rm I}_5^{(5)}|^r\right)^{1/r}\leq
C\prod_{j=1}^2\left(\sum_{|\alpha_j|=m_j}\|D^{\alpha_j}A_j\|_{\rm
BMO}\right)\||f|_r\|_{L^{n/\delta}}.
\end{equation*}
For ${\rm I}_5^{(6)}$, taking $t_1, t_2>1$ such that
$\delta/n+1/t_1+1/t_2=1$, then
\begin{align*}
\hskip -4pc \left(\sum_{i=1}^\infty|{\rm I}_5^{(6)}|^r\right)^{1/r}
&\leq C\sum_{|\alpha_1|=m_1,|\alpha_2|=m_2}\int_{R^n}\left|\frac{(x-y)^{\alpha_1+\alpha_2}K(x,y)}{|x-y|^m} \right.\\[.4pc]
&\quad\, \left. -\frac{(x_0-y)^{\alpha_1+\alpha_2}K(x_0,y)}
{|x_0-y|^m}\right| |D^{\alpha_1}\tilde A_1(y)||D^{\alpha_2}\tilde A_2(y)||h(y)|_r{\rm d}y \\[.4pc]
\hskip -4pc &\leq C\sum_{|\alpha_1|=m_1,|\alpha_2|=m_2}\sum_{k=1}^\infty(2^{-k}+2^{-\varepsilon k})\||f|_r\|_{L^{n/\delta}}  \\[.4pc]
\hskip -4pc &\quad\, \times\left(\frac{1}{|2^k\tilde Q|}\int_{2^k\tilde Q}|D^{\alpha_1}\tilde A_1(y)|^{t_1}{\rm d}y\right)^{1/t_1}\\[.4pc]
\hskip -4pc &\quad\, \times \left(\frac{1}{|2^k\tilde Q|}\int_{2^k\tilde Q}|D^{\alpha_2}\tilde A_2(y)|^{t_2}{\rm d}y\right)^{1/t_2}  \\[.4pc]
\hskip -4pc &\leq C\prod_{j=1}^2\left(\sum_{|\alpha_j|=m_j}\|D^{\alpha_j}A_j\|_{\rm BMO}\right)\||f|_r\|_{L^{n/\delta}}.
\end{align*}
Thus
\begin{equation*}
|{\rm I}_5|\leq
C\prod_{j=1}^2\left(\sum_{|\alpha_j|=m_j}\|D^{\alpha_j}A_j\|_{\rm
BMO}\right)\||f|_r\|_{L^{n/\delta}}.
\end{equation*}}
\end{pot}

\begin{pot}{\rm It suffices to prove that there exists a constant $C_Q$
such that
\begin{equation*}
\frac{1}{|Q|}\int_Q||T_\delta^A(f)(x)|_r-C_Q|{\rm d}x \leq
C\||f|_r\|_{B_p^{\delta}}
\end{equation*}
holds for any cube $Q=Q(0, d)$ with $d>1$. Without loss of generality,
we may assume $l=2$. Fix a cube $Q=Q(0, d)$ with $d>1$. Let $\tilde Q$
and $\tilde A_j(x)$ be the same as the proof of (a). Write, for
$f=g+h=\{g_i\}+\{h_i\}$ with $g_i=f_i\chi_{\tilde Q}$ and
$h_i=f_i\chi_{R^n\setminus\tilde Q}$,
\begin{align*}
&\frac{1}{|Q|}\int_Q||T_\delta^A(f)(x)|_r-|T_\delta^{\tilde A}(h)(0)|_r|{\rm d}x\\[.4pc]
&\quad\,\leq\frac{1}{|Q|} \int_Q\left(\sum_{i=1}^\infty |T_A(f_i)(x)-T_{\tilde A}(h_i)(0)|^r\right)^{1/r}{\rm d}x \\[.4pc]
&\quad\,\leq \frac{1}{|Q|}\int_Q\left(\sum_{i=1}^\infty\left|\int_{R^n}\frac{\prod_{j=1}^2R_{m_j}(\tilde A_j; x, y)}{|x-y|^m}K(x, y)g_i(y){\rm d}y\right|^r\right)^{1/r}{\rm d}x
\end{align*}
\begin{align*}
&\qquad\,+\frac{1}{|Q|}\int_Q\left(\sum_{i=1}^\infty\left|\sum_{|\alpha_1|=m_1}\frac{1}{\alpha_1!} \right.\right.\\[.3pc]
&\qquad\ \!\!\left. \left. \times \int_{R^n}\frac{R_{m_2}(\tilde A_2; x, y)
(x-y)^{\alpha_1}}{|x-y|^m}D^{\alpha_1}\tilde A_1(y)K(x,y)g_i(y){\rm d}y\right|^r\right)^{1/r}{\rm d}x \\[.3pc]
&\qquad\, +\frac{1}{|Q|}\int_Q\left(\sum_{i=1}^\infty\left|\sum_{|\alpha_2|=m_2}\frac{1}{\alpha_2!}\right.\right.\\[.3pc]
&\qquad\ \!\!\left. \left. \times\int_{R^n}\frac{R_{m_1}(\tilde A_1; x, y)
(x-y)^{\alpha_2}}{|x-y|^m}D^{\alpha_2}\tilde A_2(y)K(x,y)g_i(y){\rm d}y\right|^r\right)^{1/r}{\rm d}x \\[.3pc]
&\qquad\, +\frac{1}{|Q|}\int_Q\left(\sum_{i=1}^\infty\left|\sum_{|\alpha_1|=m_1,\ |\alpha_2|=m_2}\frac{1}{\alpha_1!\alpha_2!}\right.\right.\\[.3pc]
&\qquad\ \!\!\left. \left. \times \int_{R^n}\frac{(x-y)^{\alpha_1+\alpha_2}D^{\alpha_1}\tilde A_1(y)D^{\alpha_2}\tilde A_2(y)}{|x-y|^m}K(x,y)g_i(y){\rm d}y\right|^r\right)^{1/r}{\rm d}x \\[.3pc]
&\qquad\,+\frac{1}{|Q|}\int_Q\left(\sum_{i=1}^\infty \left|T_\delta^{\tilde A}(h_i)(x)-T_\delta^{\tilde A}(h_i)(0)\right|^r\right)^{1/r}{\rm d}x \\[.3pc]
&\quad\,:= {\rm J}_1+{\rm J}_2+{\rm J}_3+{\rm J}_4+{\rm J}_5.
\end{align*}
Similar to the proof of (a), we get, for $1/t=1/s-\delta/n$, $1<s<p$, $1<t_1,t_2<\infty$ and $1/t_1+1/t_2+s/p=1$,
\begin{align*}
{\rm J}_1&\leq C\prod_{j=1}^2\left(\sum_{|\alpha_j|=m_j}\|D^{\alpha_j}A_j\|_{\rm BMO}\right)\frac{1}{|Q|}\int_Q |T_\delta(g)(x)|_r{\rm d}x \\[.3pc]
&\leq C\prod_{j=1}^2\left(\sum_{|\alpha_j|=m_j}\|D^{\alpha_j}A_j\|_{\rm BMO}\right)\left(\frac{1}{|Q|}\int_Q|T_\delta(g)(x)|_r^q{\rm d}x\right)^{1/q} \\[.3pc]
&\leq C\prod_{j=1}^2\left(\sum_{|\alpha_j|=m_j}\|D^{\alpha_j}A_j\|_{\rm BMO}\right)d^{-n(1/p-\delta/n)}\||f|_r\chi_{\tilde Q}\|_{L^p}  \\[.3pc]
&\leq C\prod_{j=1}^2\left(\sum_{|\alpha_j|=m_j}\|D^{\alpha_j}A_j\|_{\rm BMO}\right)\||f|_r\|_{B_p^{\delta}},\\[.3pc]
{\rm J}_2 &\leq C\sum_{|\alpha_2|=m_2}\|D^{\alpha_2}A_2\|_{\rm BMO}\sum_{|\alpha_1|=m_1}\frac{1}{|Q|}\int_Q |T_\delta(D^{\alpha_1}\tilde A_1g)(x)|_r{\rm d}x \\[.3pc]
&\leq C\sum_{|\alpha_2|=m_2}\|D^{\alpha_2}A_2\|_{\rm BMO}\sum_{|\alpha_1|=m_1}\left(\frac{1}{|Q|}\int_{R^n}|T_\delta(D^{\alpha_1}\tilde A_1g)(x)|_r^t{\rm d}x\right)^{1/t}
\end{align*}
\begin{align*}
&\leq C\sum_{|\alpha_2|=m_2}\|D^{\alpha_2}A_2\|_{\rm BMO}\sum_{|\alpha_1|=m_1}|Q|^{-1/t}\||D^{\alpha_1}\tilde A_1g|_r\|_{L^s} \\[.4pc]
&\leq C\sum_{|\alpha_2|=m_2}\|D^{\alpha_2}A_2\|_{\rm BMO}\\[.4pc]
&\quad\, \times \sum_{|\alpha_1|=m_1} \left(\frac{1}{|Q|} \int_{\tilde Q}
|D^{\alpha_1}\tilde A_1(y)|^{ps/(p-s)}{\rm d}y\right)^{(p-s)/(ps)}\\[.6pc]
&\quad\, \times |Q|^{\delta/n-1/p}\||f|_r\chi_{\tilde Q}\|_{L^p} \\[.4pc]
&\leq C\prod_{j=1}^2\left(\sum_{|\alpha_j|=m_j}\|D^{\alpha_j}A_j\|_{\rm BMO}\right)\||f|_r\|_{B_p^{\delta}}, \\[.4pc]
{\rm J}_3&\leq C\prod_{j=1}^2\left(\sum_{|\alpha_j|=m_j}\|D^{\alpha_j}A_j\|_{\rm BMO}\right)\||f|_r\|_{B_p^{\delta}},\\[.4pc]
{\rm J}_4&\leq C\sum_{|\alpha_1|=m_1, |\alpha_2|=m_2}\frac{1}{|Q|}\int_Q |T_\delta(D^{\alpha_1}\tilde A_1D^{\alpha_2}\tilde A_2g)(x)|_r{\rm d}x \\[.4pc]
&\leq C\sum_{|\alpha_1|=m_1, |\alpha_2|=m_2}\left(\frac{1}{|Q|}\int_{R^n}|T_\delta(D^{\alpha_1}\tilde A_1D^{\alpha_2}\tilde A_2f_1)(x)|_r^q{\rm d}x\right)^{1/q} \\[.4pc]
&\leq C\sum_{|\alpha_1|=m_1, |\alpha_2|=m_2}|Q|^{-1/s}\left(\int_{R^n}|D^{\alpha_1}\tilde A_1(x)D^{\alpha_2}\tilde A_2(x)g(x)|_r^s{\rm d}x\right)^{1/s} \\[.4pc]
&\leq C\sum_{|\alpha_1|=m_1}\left(\frac{1}{|Q|}\int_{\tilde Q}|D^{\alpha_1}\tilde A_1(x)|^{st_1}{\rm d}x\right)^{1/st_1}\\[.4pc]
&\quad\, \times \sum_{|\alpha_2|=m_2}\left(\frac{1}{|Q|}\int_{\tilde Q}|D^{\alpha_2}\tilde A_2(x)|^{st_2}{\rm d}x\right)^{1/st_2}|Q|^{\delta/n-1/p}\||f|_r\chi_{\tilde Q}\|_{L^p} \\[.4pc]
&\leq C\prod_{j=1}^2\left(\sum_{|\alpha_j|=m_j}\|D^{\alpha_j}A_j\|_{\rm BMO}\right)\||f|_r\|_{B_p^{\delta}}.
\end{align*}
For ${\rm J}_5$, we write, for $x\in Q$,
\begin{align*}
&T_\delta^{\tilde A}(h_i)(x)-T_\delta^{\tilde A}(h_i)(0)\\[.3pc]
&\quad\,=\int_{R^n}\left(\frac{K(x,y)}{|x-y|^m}-\frac{K(0,y)}{|y|^m}\right)\prod_{j=1}^2R_{m_j}(\tilde A_j; x, y)h_i(y){\rm d}y \\[.3pc]
&\qquad\,+\int_{R^n}(R_{m_1}(\tilde A_1; x, y)-R_{m_1}(\tilde A_1; 0, y))\frac{R_{m_2}(\tilde A_2; x, y)}{|y|^m}K(0, y)h_i(y){\rm d}y
\end{align*}
\begin{align*}
&\qquad\,+\int_{R^n}(R_{m_2}(\tilde A_2; x, y)-R_{m_2}(\tilde A_2; 0, y))\frac{R_{m_1}(\tilde A_1; 0, y)}{|y|^m}K(0, y)h_i(y){\rm d}y \\[.5pc]
&\qquad\,-\sum_{|\alpha_1|=m_1}\frac{1}{\alpha_1!}\int_{R^n}\left[\frac{R_{m_2}(\tilde A_2; x, y)(x-y)^{\alpha_1}}{|x-y|^m}K(x,y) \right.\\[.5pc]
&\qquad\quad\, \left. -\frac{R_{m_2}(\tilde A_2; 0, y)(-y)^{\alpha_1}}{|y|^m}K(0,y)\right]D^{\alpha_1}\tilde A_1(y)h_i(y){\rm d}y  \\[.5pc]
&\qquad\,-\sum_{|\alpha_2|=m_2}\frac{1}{\alpha_2!}\int_{R^n}\left[\frac{R_{m_1}(\tilde A_1; x, y)(x-y)^{\alpha_2}}{|x-y|^m}K(x,y)\right.\\[.5pc]
&\qquad\quad\, \left. -\frac{R_{m_1}(\tilde A_1; 0, y)(-y)^{\alpha_2}}{|y|^m}K(0,y)\right]D^{\alpha_2}\tilde A_2(y)h_i(y){\rm d}y  \\[.5pc]
&\qquad\,+\sum_{|\alpha_1|=m_1,\ |\alpha_2|=m_2}\frac{1}{\alpha_1!\alpha_2!}\int_{R^n}\left[\frac{(x-y)^{\alpha_1+\alpha_2}}{|x-y|^m}K(x,y)\right.\\[.5pc]
&\qquad\quad\, \left. -\frac{(-y)^{\alpha_1+\alpha_2}}{|y|^m}K(0,y)\right]D^{\alpha_1}\tilde A_1(y)D^{\alpha_2}\tilde A_2(y)h_i(y){\rm d}y  \\[.6pc]
&\quad\,= {\rm J}_5^{(1)}+{\rm J}_5^{(2)}+{\rm J}_5^{(3)}+{\rm J}_5^{(4)}+{\rm J}_5^{(5)}+{\rm J}_5^{(6)}.
\end{align*}
Similar to the proof of (a), we get, for $1<t_1,t_2<\infty$ and
$1/t_1+1/t_2+1/p=1$,
\begin{align*}
\hskip -4pc \left(\sum_{i=1}^\infty|{\rm J}_5^{(1)}|^r\right)^{1/r}&\leq
C\int_{R^n}\left(\frac{|x|}{|y|^{m+n+1-
\delta}}+\frac{|x|^\varepsilon}{|y|^{m+n+\varepsilon-\delta}}\right)
\prod_{j=1}^2|R_{m_j}(\tilde A_j; x, y)||h(y)|_r{\rm d}y \\[.4pc]
\hskip -4pc &\leq C\prod_{j=1}^2\left(\sum_{|\alpha_j|=m_j}\|D^{\alpha_j}A_j\|_{\rm BMO}\right)\\[.4pc]
\hskip -4pc &\quad\, \times \sum_{k=0}^\infty\int_{2^{k+1}\tilde Q\setminus2^k\tilde Q}k^2\left(\frac{|x|}
{|y|^{n+1-\delta}}+\frac{|x|^\varepsilon}{|y|^{n+\varepsilon-\delta}}\right)|f(y)|_r{\rm d}y    \\[.4pc]
\hskip -4pc &\leq C\prod_{j=1}^2\left(\sum_{|\alpha_j|=m_j}\|D^{\alpha_j}A_j\|_{\rm BMO}\right)\\[.4pc]
\hskip -4pc &\quad\, \times \sum_{k=1}^\infty k^2(2^{-k}+2^{-\varepsilon k})(2^k d)^{-n(1/p-\delta/n)}\||f|_r\chi_{2^k\tilde Q}\|_{L^p} \\[.4pc]
\hskip -4pc &\leq C\prod_{j=1}^2\left(\sum_{|\alpha_j|=m_j}\|D^{\alpha_j}A_j\|_{\rm BMO}\right)\||f|_r\|_{B_p^{\delta}},
\end{align*}
\begin{align*}
\hskip -4pc \left(\sum_{i=1}^\infty|{\rm J}_5^{(2)}|^r\right)^{1/r} &\leq C\prod_{j=1}^2\left(\sum_{|\alpha_j|=m_j}\|D^{\alpha_j}A_j\|_{\rm BMO}\right)\sum_{k=0}^\infty
\int_{2^{k+1} \tilde Q\setminus2^k\tilde Q}k\frac{|x|}{|y|^{n+1-\delta}}|f(y)|_r{\rm d}y  \\[.5pc]
\hskip -4pc &\leq C\prod_{j=1}^2\left(\sum_{|\alpha_j|=m_j}\|D^{\alpha_j}A_j\|_{\rm BMO}\right)\||f|_r\|_{B_p^{\delta}},\\[.5pc]
\hskip -4pc \left(\sum_{i=1}^\infty|{\rm J}_5^{(3)}|^r\right)^{1/r}&\leq C\prod_{j=1}^2\left(\sum_{|\alpha_j|=m_j}\|D^{\alpha_j}A_j\|_{\rm BMO}\right)\||f|_r\|_{B_p^{\delta}},\\[.5pc]
\hskip -4pc \left(\sum_{i=1}^\infty|{\rm J}_5^{(4)}|^r\right)^{1/r}&\leq C\sum_{|\alpha_1|=m_1}\int_{R^n}\left|\frac{(x-y)^{\alpha_1}K(x,y)}{|x-y|^m}-\frac{(-y)^{\alpha_1}K(0, y)}
{|y|^m}\right|\\[.5pc]
\hskip -4pc &\quad\, \times |R_{m_2}(\tilde A_2; x, y)||D^{\alpha_1}\tilde A_1(y)||h(y)|_r{\rm d}y  \\[.5pc]
\hskip -4pc &+C\sum_{|\alpha_1|=m_1}\int_{R^n}|R_{m_2}(\tilde A_2; x, y)-R_{m_2}(\tilde A_2; 0, y)|\\[.5pc]
\hskip -4pc &\quad\, \times \frac{|(-y)^{\alpha_1}K(0, y)|}
{|y|^m}|D^{\alpha_1}\tilde A_1(y)||h(y)|_r{\rm d}y  \\[.5pc]
\hskip -4pc &\leq C\sum_{|\alpha_2|=m_2}\|D^{\alpha_2}A_2\|_{\rm BMO}\\[.6pc]
\hskip -4pc &\quad\, \times \sum_{k=1}^\infty k(2^{-k}+2^{-\varepsilon k})(2^k d)^{-n(1/p-\delta/n)}\||f|_r\chi_{2^k\tilde Q}\|_{L^p} \\[.5pc]
\hskip -4pc &\times\sum_{|\alpha_1|=m_1}\left(\frac{1}{|2^k\tilde Q|}\int_{2^k\tilde Q}|D^{\alpha_1}\tilde A_1(y)|^{p'}{\rm d}y\right)^{1/p'} \\[.5pc]
\hskip -4pc &\leq C\prod_{j=1}^2\left(\sum_{|\alpha_j|=m_j}\|D^{\alpha_j}A_j\|_{\rm BMO}\right)\||f|_r\|_{B_p^{\delta}},\\[.5pc]
\hskip -4pc \left(\sum_{i=1}^\infty|{\rm J}_5^{(5)}|^r\right)^{1/r}&\leq C\prod_{j=1}^2\left(\sum_{|\alpha_j|=m_j}\|D^{\alpha_j}A_j\|_{\rm BMO}\right)\||f|_r\|_{B_p^{\delta}},\\[.5pc]
\hskip -4pc \left(\sum_{i=1}^\infty|{\rm J}_5^{(6)}|^r\right)^{1/r}&\leq C\sum_{k=1}^\infty(2^{-k}+2^{-\varepsilon k})
(2^k d)^{-n(1/p-\delta/n)}\||f|_r\chi_{2^k\tilde Q}\|_{L^p} \\[.5pc]
\hskip -4pc &\quad\,\times \sum_{|\alpha_1|=m_1}\left(\frac{1}{|2^k\tilde Q|}\int_{2^k\tilde Q}|D^{\alpha_1}\tilde A_1(y)|^{t_1}{\rm d}y\right)^{1/t_1}
\end{align*}
\begin{align*}
\hskip -4pc &\quad\,\times \sum_{|\alpha_2|=m_2}\left(\frac{1}{|2^k\tilde Q|}\int_{2^k\tilde Q}|D^{\alpha_2}\tilde A_2(y)|^{t_2}{\rm d}y\right)^{1/t_2}    \\[.4pc]
\hskip -4pc &\leq C\prod_{j=1}^2\left(\sum_{|\alpha_j|=m_j}\|D^{\alpha_j}A_j\|_{\rm BMO}\right)\||f|_r\|_{B_p^{\delta}}.
\end{align*}
Thus
\begin{equation*}
{\rm J}_5\leq C\prod_{j=1}^2 \left(\sum_{|\alpha_j|=m_j}\| D^{\alpha_j}
A_j\|_{\rm BMO}\right)\||f|_r\|_{B_p^{\delta}}.
\end{equation*}
This completes the proof of Theorem~1.}
\end{pot}

\begin{pott}{\rm It suffices to show the conclusion for the case
$\hbox{BMO}(w)=\hbox{BMO}(R^n)$ by Lemma~2, that is, it is only to prove
that there exists a constant $C_Q$ such that
\begin{equation*}
\frac{1}{|Q|}\int_Q||T^A(f)(x)|_r-C_Q|{\rm d}x \leq
C\||f|_r\|_{L^\infty(w)}
\end{equation*}
holds for any cube $Q$. Without loss of generality, we may assume $l=2$.
Fix a cube $Q=Q(x_0, d)$. Let $\tilde Q$ and $\tilde A_j(x)$ be the same
as the proof of Theorem~1. Write, for $f=g+h=\{g_i\}+\{h_i\}$ with
$g_i=f_i\chi_{\tilde Q}$ and $h_i=f_i\chi_{R^n\setminus\tilde Q}$,
\begin{align*}
&\frac{1}{|Q|}\int_Q||T^A(f)(x)|_r-|T^{\tilde A}(h)(x_0)|_r|{\rm d}x\\[.3pc]
&\quad\, \leq\frac{1}{|Q|}
\int_Q\left(\sum_{i=1}^\infty|T_A(f_i)(x)-T_{\tilde A}(h_i)(x_0)|^r\right)^{1/r}{\rm d}x \\[.3pc]
&\quad\, \leq \frac{1}{|Q|}\int_Q\left(\sum_{i=1}^\infty\left|\int_{R^n}\frac{\prod_{j=1}^2R_{m_j}(\tilde A_j; x, y)}{|x-y|^m}K(x, y)g_i(y){\rm d}y\right|^r\right)^{1/r}{\rm d}x \\[.3pc]
&\qquad\, +\frac{1}{|Q|}\int_Q\left(\sum_{i=1}^\infty\left|\sum_{|\alpha_1|=m_1}\frac{1}{\alpha_1!} \right.\right.\\[.3pc]
&\qquad\quad\, \left.\left. \times \int_{R^n}\frac{R_{m_2}(\tilde A_2; x, y)
(x-y)^{\alpha_1}}{|x-y|^m}D^{\alpha_1}\tilde A_1(y)K(x,y)g_i(y){\rm d}y\right|^r\right)^{1/r}{\rm d}x\\[.3pc]
&\qquad\, +\frac{1}{|Q|}\int_Q\left(\sum_{i=1}^\infty\left|\sum_{|\alpha_2|=m_2}\frac{1}{\alpha_2!}\right.\right.\\[.3pc]
&\qquad\quad\, \left.\left. \times \int_{R^n}\frac{R_{m_1}(\tilde A_1; x, y)
(x-y)^{\alpha_2}}{|x-y|^m}D^{\alpha_2}\tilde A_2(y)K(x,y)g_i(y){\rm d}y\right|^r\right)^{1/r}{\rm d}x 
\end{align*}
\begin{align*}
&\qquad\, +\frac{1}{|Q|}\int_Q\left(\sum_{i=1}^\infty\left|\sum_{|\alpha_1|=m_1,\ |\alpha_2|=m_2}\frac{1}{\alpha_1!\alpha_2!}\right.\right.\\[.3pc]
&\qquad\quad\, \left.\left. \times \int_{R^n}\frac{(x-y)^{\alpha_1+\alpha_2}D^{\alpha_1}\tilde A_1(y)D^{\alpha_2}\tilde A_2(y)}{|x-y|^m}K(x,y)g_i(y){\rm d}y\right|^r\right)^{1/r}{\rm d}x \\[.3pc]
&\qquad\, +\frac{1}{|Q|}\int_Q\left(\sum_{i=1}^\infty |T_\delta^{\tilde A}(h_i)(x)-T_\delta^{\tilde A}(h_i)(x_0)|^r\right)^{1/r}{\rm d}x \\[.3pc]
&:= {\rm L}_1+{\rm L}_2+{\rm L}_3+{\rm L}_4+{\rm L}_5.
\end{align*}
By the $L^s$-boundedness of $|T|_r$ for $1<s\leq \infty$ and using the
same argument as in the proof of Theorem~1, we get, for $1/t_1+1/t_2=1$,
\begin{align*}
{\rm L}_1 &\leq C\prod_{j=1}^2\left(\sum_{|\alpha_j|=m_j}\|D^{\alpha_j}A_j\|_{\rm BMO}\right)\frac{1}{|Q|}\int_Q |T(g)(x)|_r{\rm d}x \\[.4pc]
&\leq C\prod_{j=1}^2\left(\sum_{|\alpha_j|=m_j}\|D^{\alpha_j}A_j\|_{\rm BMO}\right)\||T(g)|_r\|_{L^\infty} \\[.4pc]
&\leq C\prod_{j=1}^2\left(\sum_{|\alpha_j|=m_j}\|D^{\alpha_j}A_j\|_{\rm BMO}\right)\||f|_r\|_{L^\infty(w)}, \\[.4pc]
{\rm L}_2 &\leq C\sum_{|\alpha_2|=m_2}\|D^{\alpha_2}A_2\|_{\rm BMO}\sum_{|\alpha_1|=m_1}\frac{1}{|Q|}\int_Q |T(D^{\alpha_1}\tilde A_1g)(x)|_r{\rm d}x \\[.4pc]
&\leq C\sum_{|\alpha_2|=m_2}\|D^{\alpha_2}A_2\|_{\rm BMO}\sum_{|\alpha_1|=m_1}\left(\frac{1}{|Q|}\int_{R^n}|T(D^{\alpha_1}\tilde A_1g)(x)|_r^s{\rm d}x\right)^{1/s} \\[.4pc]
&\leq C\sum_{|\alpha_2|=m_2}\|D^{\alpha_2}A_2\|_{\rm BMO}\sum_{|\alpha_1|=m_1}\left(\frac{1}{|Q|}\int_{R^n}|D^{\alpha_1}\tilde A_1(x)g(x)|_r^s{\rm d}x\right)^{1/s} \\[.4pc]
&\leq C\sum_{|\alpha_2|=m_2}\|D^{\alpha_2}A_2\|_{\rm BMO}\\[.4pc]
&\quad\, \times \sum_{|\alpha_1|=m_1}\left(\frac{1}{|Q|}\int_{\tilde Q}|D^{\alpha_1}A_1(x)
-(D^{\alpha_1}A_1)_{\tilde Q}|^s{\rm d}x\right)^{1/s}\||f|_r\|_{L^\infty} \\[.4pc]
&\leq C\prod_{j=1}^2\left(\sum_{|\alpha_j|=m_j}\|D^{\alpha_j}A_j\|_{\rm BMO}\right)\||f|_r\|_{L^\infty(w)},\\[.4pc]
{\rm L}_3 &\leq C\prod_{j=1}^2\left(\sum_{|\alpha_j|=m_j}\|D^{\alpha_j}A_j\|_{\rm BMO}\right)\||f|_r\|_{L^\infty(w)},
\end{align*}
\begin{align*}
{\rm L}_4 &\leq C\sum_{|\alpha_1|=m_1, |\alpha_2|=m_2}\frac{1}{|Q|}\int_Q |T(D^{\alpha_1}\tilde A_1D^{\alpha_2}\tilde A_2g)(x)|_r{\rm d}x \\[.4pc]
&\leq C\sum_{|\alpha_1|=m_1, |\alpha_2|=m_2}\left(\frac{1}{|Q|}\int_{R^n}|T(D^{\alpha_1}\tilde A_1D^{\alpha_2}\tilde A_2g)(x)|_r^s{\rm d}x\right)^{1/s} \\[.4pc]
&\leq C\sum_{|\alpha_1|=m_1, |\alpha_2|=m_2}w(Q)^{-1/s}\left(\int_{R^n}|D^{\alpha_1}\tilde A_1(x)D^{\alpha_2}\tilde A_2(x)g(x)|_r^s{\rm d}x\right)^{1/s} \\[.4pc]
&\leq C\sum_{|\alpha_1|=m_1, |\alpha_2|=m_2}\left(\frac{1}{|Q|}\int_{\tilde Q}|D^{\alpha_1}\tilde A_1(x)|^{st_1}{\rm d}x\right)^{1/st_1}\\[.4pc]
&\quad\, \times \left(\frac{1}{|Q|}\int_{\tilde Q}|D^{\alpha_2}\tilde A_2(x)|^{st_2}{\rm d}x\right)^{1/st_2}\||f|_r\|_{L^\infty} \\[.4pc]
&\leq C\prod_{j=1}^2\left(\sum_{|\alpha_j|=m_j}\|D^{\alpha_j}A_j\|_{\rm BMO}\right)\||f|_r\|_{L^\infty(w)},\\[.4pc]
{\rm L}_5&\leq C\prod_{j=1}^2\left(\sum_{|\alpha_j|=m_j}\|D^{\alpha_j}A_j\|_{\rm BMO}\right)\||f|_r\|_{L^\infty(w)}.
\end{align*}}
\end{pott}

\begin{pott}{\rm It suffices to prove that there exists a constant $C_Q$
such that
\begin{equation*}
\frac{1}{w(Q)}\int_Q||T^A(f)(x)|_r-C_Q|w(x){\rm d}x \leq
C\||f|_r\|_{B_p(w)}
\end{equation*}
holds for any cube $Q=Q(0, d)$ with $d>1$. Without loss of generality,
we may assume $l=2$. Fix a cube $Q=Q(0, d)$ with $d>1$. Let $\tilde Q$
and $\tilde A_j(x)$ be the same as the proof of Theorem 1. Write, for
$f=g+h=\{g_i\}+\{h_i\}$ with $g_i= f_i\chi_{\tilde Q}$ and
$h_i=f_i\chi_{R^n\setminus\tilde Q}$,
\begin{align*}
\hskip -4pc &\frac{1}{w(Q)}\int_Q||T^A(f)(x)|_r-|T^{\tilde A}(h)(0)|_r|w(x){\rm d}x\\[.4pc]
\hskip -4pc &\quad\,\leq\frac{1}{w(Q)}\int_Q\left(\sum_{i=1}^\infty
|T_A(f_i)(x)-T_{\tilde A}(h_i)(0)|^r\right)^{1/r}w(x){\rm d}x \\[.4pc]
\hskip -4pc &\quad\,\leq \frac{1}{w(Q)}\int_Q\left(\sum_{i=1}^\infty\left|\int_{R^n}\frac{\prod_{j=1}^2R_{m_j}(\tilde A_j; x, y)}{|x-y|^m}K(x, y)g_i(y){\rm d}y\right|^r\right)^{1/r}w(x){\rm d}x \\[.4pc]
\hskip -4pc &\qquad\,+\frac{1}{w(Q)}\int_Q\left(\sum_{i=1}^\infty\left|\sum_{|\alpha_1|=m_1}\frac{1}{\alpha_1!} \right.\right.\\[.4pc]
\hskip -4pc &\qquad\quad\, \times \left.\left. \int_{R^n}\frac{R_{m_2}(\tilde A_2; x, y)
(x-y)^{\alpha_1}}{|x-y|^m}D^{\alpha_1}\tilde A_1(y)K(x,y)g_i(y){\rm d}y\right|^r\right)^{1/r}w(x){\rm d}x 
\end{align*}
\begin{align*}
\hskip -4pc &\qquad\,+\frac{1}{w(Q)}\int_Q\left(\sum_{i=1}^\infty\left|\sum_{|\alpha_2|=m_2}\frac{1}{\alpha_2!}\right.\right.\\[.4pc]
\hskip -4pc &\qquad\quad\, \times \left.\left. \int_{R^n}\frac{R_{m_1}(\tilde A_1; x, y)
(x-y)^{\alpha_2}}{|x-y|^m}D^{\alpha_2}\tilde A_2(y)K(x,y)g_i(y){\rm d}y\right|^r\right)^{1/r}w(x){\rm d}x \\[.4pc]
\hskip -4pc &\qquad\,+\frac{1}{w(Q)}\int_Q\left(\sum_{i=1}^\infty\left|\sum_{|\alpha_1|=m_1,\ |\alpha_2|=m_2}\frac{1}{\alpha_1!\alpha_2!}\right.\right.\\[.4pc]
\hskip -4pc &\qquad\quad\, \times \left.\left. \int_{R^n}\frac{(x-y)^{\alpha_1+\alpha_2}D^{\alpha_1}\tilde A_1(y)D^{\alpha_2}\tilde A_2(y)}{|x-y|^m}K(x,y)g_i(y){\rm d}y\right|^r\right)^{1/r}w(x){\rm d}x \\[.4pc]
\hskip -4pc &\qquad\,+\frac{1}{w(Q)}\int_Q\left(\sum_{i=1}^\infty |T_\delta^{\tilde A}(h_i)(x)-T_\delta^{\tilde A}(h_i)(0)|^r\right)^{1/r}w(x){\rm d}x \\[.4pc]
\hskip -4pc &\quad\,:= {\rm M}_1+{\rm M}_2+{\rm M}_3+{\rm M}_4+{\rm M}_5.
\end{align*}
Similar to the proof of Theorem 1, we get
\begin{align*}
{\rm M}_1&\leq C\prod_{j=1}^2\left(\sum_{|\alpha_j|=m_j}\|D^{\alpha_j}A_j\|_{\rm BMO}\right)\frac{1}{w(Q)}\int_Q |T(g)(x)|_rw(x){\rm d}x \\[.3pc]
&\leq C\prod_{j=1}^2\left(\sum_{|\alpha_j|=m_j}\|D^{\alpha_j}A_j\|_{\rm BMO}\right)\left(\frac{1}{w(Q)}\int_Q|T(g)(x)|_r^pw(x){\rm d}x\right)^{1/p} \\[.3pc]
&\leq C\prod_{j=1}^2\left(\sum_{|\alpha_j|=m_j}\|D^{\alpha_j}A_j\|_{\rm BMO}\right)w(\tilde Q)^{-1/p}\||f|_r\chi_{\tilde Q}\|_{L^p(w)} \\[.3pc]
&\leq C\prod_{j=1}^2\left(\sum_{|\alpha_j|=m_j}\|D^{\alpha_j}A_j\|_{\rm BMO}\right)\||f|_r\|_{B_p(w)}.
\end{align*}
For ${\rm M}_2$, since $w\in A_1$, $w$ satisfies the reverse of H\"older's
inequality:
\begin{equation*}
\left(\frac{1}{|Q|}\int_Q w(x)^q{\rm d}x\right)^{1/q}\leq
\frac{C}{|Q|}\int_Q w(x){\rm d}x
\end{equation*}
for all cube $Q$ and some $1<q<\infty$ \cite{10,16}. Thus, taking $s,
t>1$ such that $st<p$ and $q=(pt-st)/(p-st)$, then
\begin{align*}
{\rm M}_2 &\leq C\sum_{|\beta_2|=m_2}\|D^{\alpha_2}A_2\|_{\rm BMO}\sum_{|\alpha_1|=m_1}\frac{1}{w(Q)}\int_Q |T(D^{\alpha_1}\tilde A_1g)(x)|_rw(x){\rm d}x \\[.3pc]
&\leq C\sum_{|\alpha_2|=m_2}\|D^{\alpha_2}A_2\|_{\rm BMO}
\end{align*}
\begin{align*}
&\quad\, \times \sum_{|\alpha_1|=m_1}\left(\frac{1}{w(Q)}\int_{R^n}|T(D^{\alpha_1}\tilde A_1g)(x)|_r^sw(x){\rm d}x\right)^{1/s}\\[.3pc]
&\leq C\sum_{|\alpha_2|=m_2}\|D^{\alpha_2}A_2\|_{\rm BMO}\sum_{|\alpha_1|=m_1}w(Q)^{-1/s}\|D^{\alpha_1}\tilde A_1|g|_r\|_{L^s(w)} \\[.3pc]
&\leq C\sum_{|\alpha_2|=m_2}\|D^{\alpha_2}A_2\|_{\rm BMO}w(Q)^{-1/s}\sum_{|\alpha_1|=m_1}\left(\int_{\tilde Q}|D^\alpha\tilde A_1(y)|^{st'}{\rm d}y\right)^{1/st'}\\[.3pc]
&\quad\, \times \left(\int_{\tilde Q}|f(x)|_r^{st}w(x)^t{\rm d}x\right)^{1/st} \\[.3pc]
&\leq C\prod_{j=1}^2\left(\sum_{|\alpha_j|=m_j}\|D^{\alpha_j}A_j\|_{\rm BMO}\right)|Q|^{1/st'}w(Q)^{-1/s}\\[.3pc]
&\quad\, \times \left(\int_{\tilde Q}|f(x)|_r^pw(x){\rm d}x\right)^{1/p}\left(\int_{\tilde Q}w(x)^q{\rm d}x\right)^{(p-s)/pqs}  \\[.3pc]
&\leq C\prod_{j=1}^2\left(\sum_{|\alpha_j|=m_j}\|D^{\alpha_j}A_j\|_{\rm BMO}\right)w(\tilde Q)^{-1/p}\||f|_r\chi_{\tilde Q}\|_{L^p(w)} \\[.3pc]
&\leq C\prod_{j=1}^2\left(\sum_{|\alpha_j|=m_j}\|D^{\alpha_j}A_j\|_{\rm BMO}\right)\||f|_r\|_{B_p(w)}, \\[.3pc]
{\rm M}_3 &\leq C\prod_{j=1}^2\left(\sum_{|\alpha_j|=m_j}\|D^{\alpha_j}A_j\|_{\rm BMO}\right)\||f|_r\|_{B_p(w)}.
\end{align*}
For ${\rm M}_4$, taking $s, t_1, t_2, t_3>1$ such that $1/t_1+1/t_2+1/t_3=1$,
$st_3<p$ and $q=(pt_3-st_3)/(p-st_3)$, then, by the reverse of
H\"older's inequality,
\begin{align*}
\hskip -4pc {\rm M}_4 &\leq C\sum_{|\alpha_1|=m_1, |\alpha_2|=m_2}\frac{1}{w(Q)}\int_Q |T(D^{\alpha_1}\tilde A_1D^{\alpha_2}\tilde A_2g)(x)|_rw(x){\rm d}x \\[.3pc]
\hskip -4pc &\leq C\sum_{|\alpha_1|=m_1, |\alpha_2|=m_2}\left(\frac{1}{w(Q)}\int_{R^n}|T(D^{\alpha_1}\tilde A_1D^{\alpha_2}\tilde A_2g)(x)|_r^sw(x){\rm d}x\right)^{1/s} \\[.3pc]
\hskip -4pc &\leq C\sum_{|\alpha_1|=m_1, |\alpha_2|=m_2}w(Q)^{-1/s}\left(\int_{R^n}|D^{\alpha_1}\tilde A_1(x)D^{\alpha_2}\tilde A_2(x)g(x)|_r^sw(x){\rm d}x\right)^{1/s} \\[.3pc]
\hskip -4pc &\leq C\sum_{|\alpha_1|=m_1}\left(\int_{\tilde Q}|D^{\alpha_1}\tilde A_1(x)|^{st_1}{\rm d}x\right)^{1/st_1}\\[.3pc]
\hskip -4pc &\quad\times \sum_{|\alpha_2|=m_2}\left(\int_{\tilde Q}|D^{\alpha_2}\tilde A_2(x)|^{st_2}{\rm d}x\right)^{1/st_2} 
w(Q)^{-1/s} \left(\int_{\tilde Q}|f(x)|_r^{st_3}w(x)^{t_3}{\rm d}x\right)^{1/st_3} 
\end{align*}
\begin{align*}
\hskip -4pc &\leq C\prod_{j=1}^2\left(\sum_{|\alpha_j|=m_j}\|D^{\alpha_j}A_j\|_{\rm BMO}\right)w(\tilde Q)^{-1/p}\||f|_r\chi_{\tilde Q}\|_{L^p(w)} \\[.3pc]
\hskip -4pc &\leq C\prod_{j=1}^2\left(\sum_{|\alpha_j|=m_j}\|D^{\alpha_j}A_j\|_{\rm BMO}\right)\||f|_r\|_{B_p(w)}.
\end{align*}
For ${\rm M}_5$, we write, for $x\in Q$,
\begin{align*}
&T^{\tilde A}(h_i)(x)-T^{\tilde A}(h_i)(0)\\[.4pc]
&\quad\,=\int_{R^n}\left(\frac{K(x,y)}{|x-y|^m}-\frac{K(0,y)}{|y|^m}\right)\prod_{j=1}^2R_{m_j}(\tilde A_j; x, y)h_i(y){\rm d}y \\[.5pc]
&\qquad\,+\int_{R^n}(R_{m_1}(\tilde A_1; x, y)-R_{m_1}(\tilde A_1; 0, y))\frac{R_{m_2}(\tilde A_2; x, y)}{|y|^m}K(0, y)h_i(y){\rm d}y \\[.5pc]
&\qquad\,+\int_{R^n}(R_{m_2}(\tilde A_2; x, y)-R_{m_2}(\tilde A_2; 0, y))\frac{R_{m_1}(\tilde A_1; 0, y)}{|y|^m}K(0, y)h_i(y){\rm d}y \\[.5pc]
&\qquad\,-\sum_{|\alpha_1|=m_1}\frac{1}{\alpha_1!}\int_{R^n}\left[\frac{R_{m_2}(\tilde A_2; x, y)(x-y)^{\alpha_1}}{|x-y|^m}K(x,y) \right.\\[.5pc]
&\qquad\quad\,\left. -\frac{R_{m_2}(\tilde A_2; 0, y)(-y)^{\alpha_1}}{|y|^m}K(0,y)\right]D^{\alpha_1}\tilde A_1(y)h_i(y){\rm d}y  \\[.5pc]
&\qquad\,-\sum_{|\alpha_2|=m_2}\frac{1}{\alpha_2!}\int_{R^n}\left[\frac{R_{m_1}(\tilde A_1; x, y)(x-y)^{\alpha_2}}{|x-y|^m}K(x,y)\right.\\[.5pc]
&\qquad\quad\,\left. -\frac{R_{m_1}(\tilde A_1; 0, y)(-y)^{\alpha_2}}{|y|^m}K(0,y)\right]D^{\alpha_2}\tilde A_2(y)h_i(y){\rm d}y  \\[.5pc]
&\qquad\,+\sum_{|\alpha_1|=m_1,\ |\alpha_2|=m_2}\frac{1}{\alpha_1!\alpha_2!}\int_{R^n}\left[\frac{(x-y)^{\alpha_1+\alpha_2}}{|x-y|^m}K(x,y)\right.\\[.5pc]
&\qquad\quad\,\left. -\frac{(-y)^{\alpha_1+\alpha_2}}{|y|^m}K(0,y)\right]D^{\alpha_1}\tilde A_1(y)D^{\alpha_2}\tilde A_2(y)h_i(y){\rm d}y  \\[.5pc]
&\quad\,= {\rm M}_5^{(1)}+{\rm M}_5^{(2)}+{\rm M}_5^{(3)}+{\rm M}_5^{(4)}+{\rm M}_5^{(5)}+{\rm M}_5^{(6)}.
\end{align*}
Similar to the proof of Theorem~1 and notice that $w\in A_1\subset A_p$,
we get
\begin{align*}
\hskip -4pc\left(\sum_{i=1}^\infty|{\rm M}_5^{(1)}|^r\right)^{1/r}&\leq C\int_{R^n}\left(\frac{|x|}{|y|^{m+n+1}}+\frac{|x|^\varepsilon}
{|y|^{m+n+\varepsilon}}\right)\prod_{j=1}^2|R_{m_j}(\tilde A_j; x, y)||h(y)|_r{\rm d}y 
\end{align*}
\begin{align*}
\hskip -4pc&\leq C\prod_{j=1}^2\left(\sum_{|\alpha_j|=m_j}\|D^{\alpha_j}A_j\|_{\rm BMO}\right)\\[.4pc]
\hskip -4pc&\quad\, \times \sum_{k=0}^\infty\int_{2^{k+1}\tilde Q\setminus2^k\tilde Q}
k^2\left(\frac{|x|}{|y|^{n+1}}+\frac{|x|^\varepsilon}{|y|^{n+\varepsilon}}\right)|f(y)|_r{\rm d}y\\[.4pc]
\hskip -4pc&\leq C\prod_{j=1}^2\left(\sum_{|\alpha_j|=m_j}\|D^{\alpha_j}A_j\|_{\rm BMO}\right)\\[.4pc]
\hskip -4pc &\quad\, \times \sum_{k=1}^\infty k^2(2^{-k}+2^{-\varepsilon k})
w(2^k\tilde Q)^{-1/p}\left(\int_{2^k\tilde Q}|f(y)|_r^pw(y){\rm d}y\right)^{1/p}  \\[.4pc]
\hskip -4pc&\times\left(\frac{1}{|2^k\tilde Q|}\int_{2^k\tilde Q}w(y){\rm d}y\right)^{1/p}\left(\frac{1}{|2^k\tilde Q|}\int_{2^k\tilde Q}w(y)^{-1/(p-1)}{\rm d}y\right)^{(p-1)/p} \\[.4pc]
\hskip -4pc&\leq C\prod_{j=1}^2\left(\sum_{|\alpha_j|=m_j}\|D^{\alpha_j}A_j\|_{\rm BMO}\right)\||f|_r\|_{B_p(w)},\\[.4pc]
\hskip -4pc\left(\sum_{i=1}^\infty|{\rm M}_5^{(2)}|^r\right)^{1/r}&\leq C\prod_{j=1}^2\left(\sum_{|\alpha_j|=m_j}\|D^{\alpha_j}A_j\|_{\rm BMO}\right)
\sum_{k=0}^\infty\int_{2^{k+1}\tilde Q\setminus2^k\tilde Q}k\frac{|x|}{|y|^{n+1}}|f(y)|_r{\rm d}y  \\[.4pc]
\hskip -4pc&\leq C\prod_{j=1}^2\left(\sum_{|\alpha_j|=m_j}\|D^{\alpha_j}A_j\|_{\rm BMO}\right)\||f|_r\|_{B_p(w)}, \\[.4pc]
\hskip -4pc\left(\sum_{i=1}^\infty|{\rm M}_5^{(3)}|^r\right)^{1/r}&\leq C\prod_{j=1}^2\left(\sum_{|\alpha_j|=m_j}\|D^{\alpha_j}A_j\|_{\rm BMO}\right)\||f|_r\|_{B_p(w)}.
\end{align*}
For ${\rm M}_5^{(4)}$, choose $1<s<p$, notice that $w\in A_1\subset A_{p/s}$, we get
\begin{align*}
\hskip -4pc\left(\sum_{i=1}^\infty|{\rm M}_5^{(4)}|^r\right)^{1/r}&\leq C\sum_{|\alpha_2|=m_2}\|D^{\alpha_2}A_2\|_{\rm BMO}\\[.4pc]
\hskip -4pc &\quad\, \times \sum_{k=0}^\infty
\int_{2^{k+1}\tilde Q\setminus2^k\tilde Q}k\left(\frac{|x|}{|y|^{n+1}}+\frac{|x|^\varepsilon}{|y|^{n+\varepsilon}}\right)|D^{\alpha_1}\tilde A_1(y)||f(y)|_r{\rm d}y \\[.4pc]
\hskip -4pc&\leq C\sum_{|\alpha_2|=m_2}\|D^{\alpha_2}A_2\|_{\rm BMO}\sum_{k=0}^\infty\left(\frac{d}{(2^kd)^{n+1}}+\frac{d^\varepsilon}{(2^kd)^{n+\varepsilon}}\right)\\[.4pc]
\hskip -4pc&\quad\, \times \left(\int_{2^{k+1}\tilde Q}|f(y)|_r^s{\rm d}y\right)^{1/s}{\rm d}y \left(\int_{2^{k+1}\tilde Q}|D^{\alpha_1}\tilde A_1(y)|^{s'}{\rm d}y\right)^{1/s'}
\end{align*}
\begin{align*}
\hskip -4pc&\leq C\prod_{j=1}^2\left(\sum_{|\alpha_j|=m_j}\|D^{\alpha_j}A_j\|_{\rm BMO}\right)\sum_{k=1}^\infty k(2^{-k}+2^{-\varepsilon k})
w(2^k\tilde Q)^{-1/p}\\[.4pc]
\hskip -4pc&\quad\, \times\left(\int_{2^k\tilde Q}|f(y)|_r^pw(y){\rm d}y\right)^{1/p} \left(\frac{1}{|2^k\tilde Q|}\int_{2^k\tilde Q}w(y){\rm d}y\right)^{1/p}\\[.4pc]
\hskip -4pc&\quad\, \times\left(\frac{1}{|2^k\tilde Q|}\int_{2^k\tilde Q}w(y)^{-s/(p-s)}{\rm d}y\right)^{(p-s)/ps} \\[.4pc]
\hskip -4pc&\leq C\prod_{j=1}^2\left(\sum_{|\alpha_j|=m_j}\|D^{\alpha_j}A_j\|_{\rm BMO}\right)\||f|_r\|_{B_p(w)}, \\[.4pc]
\hskip -4pc\left(\sum_{i=1}^\infty|{\rm M}_5^{(5)}|^r\right)^{1/r}&\leq C\prod_{j=1}^2\left(\sum_{|\alpha_j|=m_j}\|D^{\alpha_j}A_j\|_{\rm BMO}\right)\||f|_r\|_{B_p(w)}.
\end{align*}
For ${\rm L}_5^{(6)}$, choose $1<t_1,t_2,t_3<\infty$ such that $t_3<p$ and $1/t_1+1/t_2+1/t_3=1$. Notice that $w\in A_1\subset A_{p/t_3}$, we get
\begin{align*}
\hskip -4pc\left(\sum_{i=1}^\infty|{\rm M}_5^{(6)}|^r\right)^{1/r}
&\leq C\sum_{|\alpha_1|=m_1,|\alpha_2|=m_2}\int_{R^n}\left|\frac{(x-y)^{\alpha_1+\alpha_2}K(x,y)}{|x-y|^m} -\frac{(-y)^{\alpha_1+\alpha_2}K(0, y)}
{|y|^m}\right|\\[.4pc]
&\quad\, \times |D^{\alpha_1}\tilde A_1(y)||D^{\alpha_2}\tilde A_2(y)||h(y)|_r{\rm d}y \\[.4pc]
\hskip -4pc&\leq C\sum_{k=0}^\infty\left(\frac{d}{(2^kd)^{n+1}}+\frac{d^\varepsilon}{(2^kd)^{n+\varepsilon}}\right)
\left(\int_{2^{k+1}\tilde Q}|f(y)|_r^{t_3}{\rm d}y\right)^{1/t_3}{\rm d}y \\[.4pc]
\hskip -4pc&\times\sum_{|\alpha_1|=m_1}\left(\int_{2^{k+1}\tilde Q}|D^{\alpha_1}\tilde A_1(y)|^{t_1}{\rm d}y\right)^{1/t_1}\\[.4pc]
&\quad\times \sum_{|\alpha_2|=m_2}\left(\int_{2^{k+1}\tilde Q}|D^{\alpha_2}\tilde A_2(y)|^{t_2}{\rm d}y\right)^{1/t_2} \\[.4pc]
\hskip -4pc&\leq C\prod_{j=1}^2\left(\sum_{|\alpha_j|=m_j}\|D^{\alpha_j}A_j\|_{\rm BMO}\right)\sum_{k=1}^\infty(2^{-k}+2^{-\varepsilon k})
w(2^k\tilde Q)^{-1/p}\\[.4pc]
&\quad\, \times \left(\int_{2^k\tilde Q}|f(y)|_r^pw(y){\rm d}y\right)^{1/p}
\left(\frac{1}{|2^k\tilde Q|}\int_{2^k\tilde Q}w(y){\rm d}y\right)^{1/p}\\[.4pc]
\hskip -4pc&\times\left(\frac{1}{|2^k\tilde Q|}\int_{2^k\tilde Q}w(y)^{-t_3/(p-t_3)}{\rm d}y\right)^{(p-t_3)/pt_3} \\[.4pc]
\hskip -4pc&\leq C\prod_{j=1}^2\left(\sum_{|\alpha_j|=m_j}\|D^{\alpha_j}A_j\|_{\rm BMO}\right)\||f|_r\|_{B_p(w)}.
\end{align*}
Thus
\begin{equation*}
{\rm M}_5\leq C\prod_{j=1}^2\left(\sum_{|\alpha_j|=m_j}\|D^{\alpha_j}A_j\|_{\rm
BMO}\right)\||f|_r\|_{B_p(w)}.
\end{equation*}
This completes the proof of Theorem~2.}\vspace{-.7pc}
\end{pott}

\section{Applications}

Now we shall apply the theorems of the paper to some particular
operators such as the Calder\'on--Zygmund singular integral operator and
fractional integral operator.

\begin{application}{\rm 
Calder\'on--Zygmund singular integral operator.}
\end{application}

Let $T$ be the Calder\'on--Zygmund operator \cite{7,10,16}. Then it is
easy to see that $T$ satisfies the conditions in Theorem~2. Thus the
conclusions of Theorem~2 hold for $T^A$.

\begin{application}{\rm 
Fractional integral operator with rough kernel.}
\end{application}

For $0<\delta<n$, let $T_\delta$ be the fractional integral operator
with rough kernel defined by \cite{1,9,11}
\begin{equation*}
T_\delta f(x)=\int_{R^n}\frac{\Omega(x-y)}{|x-y|^{n-\delta}}f(y){\rm d}y,
\end{equation*}
the vector-valued multilinear operator related to $T_\delta$ is defined
by
\begin{equation*}
|T_\delta^A(f)(x)|_r=\left(\sum_{i=1}^\infty
|T_\delta^A(f_i)(x)|^r\right)^{1/r},
\end{equation*}
where
\begin{equation*}
T_\delta^A(f)(x)=\int_{R^n}\frac{\prod_{j=1}^lR_{m_j+1}(A_j; x,
y)}{|x-y|^{m+n-\delta}}\Omega(x-y)f(y){\rm d}y
\end{equation*}
and $\Omega$ is homogeneous of degree zero on $R^n$,
$\int_{S^{n-1}}\Omega(x'){\rm d}\sigma(x')=0$ and $\Omega\in
{\rm Lip}_\varepsilon(S^{n-1})$ for some $0<\varepsilon\leq 1$, that is, there
exists a constant $M>0$ such that for any $x,y\in S^{n-1}$,
$|\Omega(x)-\Omega(y)|\leq M|x-y|^\varepsilon$. Then $T_\delta$
satisfies the conditions in Theorem~1. Thus the conclusions of Theorem~1
hold for $T_\delta^A$.\vspace{-.7pc}

\section*{Acknowledgements}

The author would like to express his deep gratitude to the referee for
his valuable comments and suggestions. This work is supported by the
NNSF (Grant: 10271071).

\end{document}